\documentclass[reqno,12pt]{amsart}

\usepackage{amsmath}
\usepackage{amsfonts}
\usepackage{amssymb}
\usepackage{eufrak}
\usepackage{amscd}
\usepackage{amsthm}
\usepackage{epsfig}
\usepackage{amstext}
\usepackage[all]{xy}

\theoremstyle{plain}

 \newtheorem{theorem}{Theorem}
 
 \newtheorem{question}[theorem]{Question}
 \newtheorem{prop}[theorem]{Proposition}
\newtheorem{opred}[theorem]{Definition}
 \newtheorem{cor}[theorem]{Corollary}
 \newtheorem{lemma} [theorem] {Lemma}
 \newtheorem{remark}[theorem]{Remark}

\newcommand{\pnl}{P_{n, \lambda}}

\begin{document}

 \title{On universal $C^*$-algebras generated by n projections with scalar sum}
\author{Tatiana Shulman}

\address{Mathematics Department, New Hampshire University, Durham,
New Hampshire 03824, USA}

\email{tatiana$_-$shulman@yahoo.com}



\date{\today}

\maketitle


\section*{Introduction}
We consider the relations
\begin{equation}\label{osnsootn}
 \sum_{i=1}^np_i = \lambda 1,
  \; p_i = p_i^* = p_i^2, \; i=\overline{1,  n},
 \end{equation}
where  $\lambda\in \mathbb{R}$, and their representations, that is
n-tuples $P_i, i=\overline{1, n},$ of projections on a Hilbert space
such that $\sum_{i=1}^n P_i = \lambda 1$.

Decompositions of scalar operators  on a Hilbert space into a sum of
fixed number of projections were studied in series of papers
(\cite{LinAlg}, \cite{KrugRoyt}, \cite{RabSam}, \cite{Krug},
\cite{Yu.S.}). In \cite{Yu.S.} it was completely described  the set
$\Sigma_n$ of all scalars $\lambda$ such that the operators $\lambda
1$ are sums of n projections. It was proved that for $n<4$,
$\;\Sigma_n$ is finite:
\begin{itemize} \label{maln}\item $\Sigma_1 = \{0,
 1\},$

 \item  $\Sigma_2 = \{0, 1, 2\},$

 \item  $\Sigma_3 = \{0, 1, 3/2, 2, 3\}$,

 \end{itemize}
  for  $n=4$, it is countable:
\begin{itemize}  \item   $\Sigma_4 = \{0, 1, 1+\frac{k}{k+2} (k\in \mathbb N), 2, 3-\frac{k}{k+2} (k\in \mathbb N), 3, 4\}$,
\end{itemize}
  and for $n\ge 5$, $\;\Sigma_n$ is the union of the "main"
  segment
 $[\alpha_n,\beta_n]$  and two sequences  $\Lambda_n^1$ and $\Lambda_n^2$ converging to the end points of this segment.
 Here \begin{equation}\label{KontsyOtrezka}
 \alpha_n = \frac{n - \sqrt{n^2-4n}}{2}, \;\;\;\; \beta_n = \frac{n + \sqrt{n^2-4n}}{2}
 \end{equation}
 We don't need the exact formula for the points from $\Lambda_n^1$ and
 $\Lambda_n^2$, it will be sufficient to know that
 they
 are rational. It is not difficult to see that $\lambda \in \Sigma_n$ if and only if the
relation (\ref{osnsootn}) has a representation in a separable
Hilbert space. The main interest is to describe, for each
$\lambda\in \Sigma_n$, all n-tuples of projections in a Hilbert
space that fulfill (\ref{osnsootn}) or at least to understand how
complicated this problem is. The
 degree of complexity of this task for arbitrary relation
 can be formulated in terms of belonging of the universal $C^*$-algebra of  the relation to some  less or more
 tractable classes of
 $C^*$-algebras (type I, approximately finite dimensional, nuclear  $C^*$-algebras). For $\lambda\in \Sigma_n$, let $\pnl$ be the universal
$C^*$-algebra of the relation (\ref{osnsootn}).

In \cite{Yu.S.} the authors ask for which $\lambda$ the
$C^*$-algebra $\mathcal{P}_{n, \lambda}$ belongs to the class of
type I $C^*$-algebras. They proved that if $\lambda\in
\Lambda_n^i, i=1, 2,$  then $\pnl$ is finite-dimensional and if
$\lambda\in (\alpha_n, \beta_n)$ then $\pnl$ is not of type I (for
any $n>6$). For $\lambda = \alpha_n$ and $\lambda = \beta_n$ the
question remained open. Below we'll give the negative answer to
this question. Moreover it will be shown that for these values of
$\lambda$  there  don't exist unital $\ast$-homomorphisms from
 $\mathcal{P}_{n, \lambda}$ to any type I $C^*$-algebra.

Then we show that for "most" values of $\lambda$ the $C^*$-algebras
$\pnl$ are not nuclear and even exact.

We prove that for every $\lambda$, $\pnl$ has a trace and use this
fact in the problem  of classification of these $C^*$-algebras. The
result is  that among these $C^*$-algebras there is continuum of
mutually nonisomorphic ones.

We also study what scalars can be represented as a sum of n
projections in a given unital $C^*$-algebra. For arbitrary unital
$C^*$-algebra $A$, let us denote the set of such scalars by
$\Sigma_n(A)$. We explore  $\Sigma_n(A)$ for different classes of
$C^*$-algebras (type I $C^*$-algebras, $C^*$-algebras with a
trace), in particular for UHF-algebras we completely describe this
set.

All these results are presented in Chapter 2. In Chapter 1 we
consider some questions of general theory of representations of
polynomial relations. We introduce a topology on the set of all
polynomial relations and are interested in the question when the set
of relations representable in a given class of $C^*$-algebras (such
as AF-algebras and $C^*$-algebras with a trace) is closed.  The
results will be applied in Chapter 2  to the relations
(\ref{osnsootn}).

 The author is grateful to Yu.S. Samoilenko for
stimulating questions and helpful discussions.

\section{ Representations of polynomial relations in AF-algebras and  $C^*$-algebras with a trace}

  By {\it polynomial relation} (in n variables) we  call an equation of the form
\begin{equation}\label{sootn}
f(x_1,...,x_n,x_1^{\ast},...,x_n^{\ast}) = 0
\end{equation}
where $f$ is a polynomial in $2n$ noncommuting variables, that is
 an element of  the free unital $\ast$-algebra $\mathcal{F}^{\ast}_{n}$ on
generators $x_1,...,x_n$. A {\it representation} of (\ref{sootn}) is
any n-tuple $(T_1,...,T_n)$ of operators in a Hilbert space that
fulfills the condition
\begin{equation}\label{predstsootn}
f(T_1,...,T_n,T_1^{\ast},...,T_n^{\ast}) = 0.
\end{equation}

We consider also  representations in unital $C^*$-algebras. Let us
say that the relation (\ref{sootn}) is {\it representable in a
(unital) $C^*$-algebra $A$} if there is an n-tuple of elements
$a_1,...,a_n \in A$ that satisfies the equality
\begin{equation}\label{starpredstsootn}
f(a_1,...,a_n,a_1^{\ast},...,a_n^{\ast}) = 0.
\end{equation}
(Strictly speaking, (\ref{starpredstsootn}) means that the image of
the element $f$ under the unital $\ast$-homomorphism of
$\mathcal{F}^{\ast}_{n}$ to $A$ which sends each $x_i$ to $a_i$, is
zero).

If $\|a_i\|\le C$ for all $i$, where $C>0$, then we say that this
representation is {\it $C$-bounded}.

The  universal $C^*$-algebra of the relation (\ref{sootn}) is a
(unital) $C^*$-algebra $A$  generated by elements  $a_1,...,a_n$
such that
 (\ref{starpredstsootn}) holds and for any representation
$(T_1,...,T_n)$ of (\ref{sootn}) there is a  $\ast$-representation
$\pi$ of $A$ such that $\pi(a_i) = T_i$.

A necessary and sufficient condition for the existence of the
universal $C^*$-algebra of a polynomial relation is
\begin{equation}\label{ogr} \sup \max_i\|T_i\|
<\infty,
\end{equation}
 where supremum  is taken over the set of all representations of the relation.

Representations and universal $C^*$-algebras of systems of
relations are defined in a similar way. Actually there is no
difference between the case of one relation and the case of a
finite system of relations because there is an easy way to rewrite
the system
\begin{equation}\label{nabsootn}
f_i(a_1,...,a_n,a_1^{\ast},...,a_n^{\ast}) = 0, (1\le i\le m)
\end{equation}
in the form of one relation
 \begin{equation}\label{sumsootn}
\sum_i(f_i(a_1,...,a_n,a_1^{\ast},...,a_n^{\ast}))^{\ast}f_i(a_1,...,a_n,a_1^{\ast},...,a_n^{\ast})
= 0.
\end{equation}
So we shall say, for short, about representability of polynomial
relations but mean that results concern also systems of polynomial
relations.

Now let us introduce a topology on the set of all polynomial
relations (which can be identified with $\mathcal{F}^{\ast}_{n}$) by
a system of seminorms ${\nu}_K$, $K > 0$ on $\mathcal{F}^{\ast}_{n}$
, defined by  the formula
$${\nu}_K(f) = \sup \|f(T_1,\ldots,T_n, T_1^*, \ldots, T_n^*)\|,$$
where supremum is taken over the set of all n-tuples of operators
with norms not greater than $K$.

For the rest of the paper all $C^*$-algebras and
$\ast$-homomorphisms are assumed to be unital.

\subsection{ Representations in $C^*$-algebras with a trace}

By a  {\it trace}  on a  $C^*$-algebra $A$ we mean a positive
linear functional $g$, such that
$$g(xy)=g(yx),\;g(1_A) = 1.$$

\begin{theorem}\label{pr-sled}
Let $C>0$. The set of polynomial relations in $n$ variables,
$C$-boundedly representable in  $C^*$-algebras with a trace, is
closed.
\end{theorem}
{\bf Proof.} Let $f_k\to f$ and for any $k$ the relation $f_k$ has a
$C$-bounded representation $\pi_k$ in a $C^*$-algebra $A_k$ with a
trace $\tau_k$.  Consider the $C^*$-algebra $E$ of all bounded
sequences
 $(a_k)_{k\in \mathbb{N}}$, where $ a_k\in A_k$. Let $J$
be the ideal in $E$ consisting of all sequences vanishing at
infinity and let $F = E/J$.  Set $$a_{i,k} = \pi_k(x_i),$$
$i=\overline{1, n}$, where $x_i$ are free generators of $\mathcal
F_n^*$ . Let $e_i = (a_{i,1},a_{i,2},...)\in E$ and $b_i$ be their
images under the canonical epimorphism from $E$ to $F$. Then,
setting
$$\pi(x_i) = b_i,$$ we define a representation of the relation $f$
in $F$.

It remains to prove that $F$ has a trace. Let $\xi$ be a
nontrivial ultrafilter on $\mathbb N$. We can think of it as a
character of $l^{\infty}(\mathbb N)$. Setting
$$h((a_k)_{k\in \mathbb{N}}) = \xi((\tau_k(a_k))_{k\in
\mathbb{N}}),$$ we get a state on  $E$. It is easy to see that
$h(xy)=h(yx)$, for any $x, y \in E$,  $h(1_E) = 1$ and hence $h$
is a trace on $E$. Since any ultrafilter vanishes on  $c_0$, $h$
vanishes on $J$ and hence defines a trace on $F$. \qed

\subsection{ Representations in AF-algebras}
 It is a consequence of Theorem \ref{pr-sled} that if a relation
 is a limit of a sequence of relations $C$-boundedly representable in
 finite-dimensional C*-algebras then it is representable in a C*-algebra with a trace.
 Now we obtain a partial converse of this result, for a more narrow
class of all AF-algebras.

Recall that a C*-algebra $A$ is called an {\it AF-algebra} if
there is an up-directed by inclusion net $A_{k}$ of finite
dimensional subalgebras, such that  $A = \overline{\bigcup A_k}$.
We assume that $A$ is unital and all $A_k$ contain the unit of
$A$. It is not difficult to see that each (unital) AF-algebra has
a trace.
\begin{lemma}\label{mnogochlen} Let  $0<\varepsilon<D$, $-\varepsilon< \lambda_1, ..., \lambda_N < \varepsilon$,
$N\in \mathbb N$. Then there exists a polynomial $h(t)$ such that

a) $|h(t) - t| < \varepsilon$ when $-D \le t \le D$;

b) $h(\lambda_i) = 0 $ for any  $1\le i \le N$.
\end{lemma}
{\bf Proof.} Let $\varepsilon_1 < \varepsilon$ be a number greater
that all $|\lambda_i|$, $i = \overline{1, N}$. Set $$f(t)=
\begin{cases} 0,\; t\in
(-\varepsilon_1,\varepsilon_1) \\ t -\varepsilon_1, \; t>\varepsilon_1 \\
t+\varepsilon_1, \;t < -\varepsilon_1 \end{cases}$$ Then  $|f(t)
-t|\le \varepsilon_1$ for all t. Set $h_0(t) =
(t-\lambda_1)\ldots(t-\lambda_N)$ and let $M= sup_{|t|\le D}
|h_0(t)|.$ Then $f(t)/h_0(t)$ is a continuous function and by
Weierstrass theorem there is a polynomial $q(t)$ such that
$$|q(t) - f(t)/h_0(t) | < (\varepsilon-\varepsilon_1)/ M$$ when
$-D\le t\le D$, whence
$$|q(t)h_0(t) - t| \le |q(t)h_0(t) - f(t)| + |f(t) - t| <
\varepsilon$$ when $-D\le t\le D$.  Set $h(t) = q(t)h_0(t)$. It is
easy to see that the conditions a) and b) are fulfilled. \qed

\begin{theorem}\label{lim}
Any polynomial relation  representable in  AF $C^*$-algebra is a
limit of polynomial relations $C$-boundedly representable in
finite-dimensional $C^*$-algebras, for some $C>0$.
\end{theorem}
{\bf Proof.} Let \begin{equation}\label{sootn'}
f(x_1,...,x_n,x_1^{\ast},...,x_n^{\ast}) = 0
\end{equation} be a polynomial relation representable in an AF-algebra.
This relation has the same representations as the relation
 $f(x_1,...,x_n,x_1^{\ast},...,x_n^{\ast})^*
f(x_1,...,x_n,x_1^{\ast},...,x_n^{\ast})= 0.$
 Hence, replacing, if necessary, $f(x_1,...,x_n,x_1^{\ast},...,x_n^{\ast})$ by
 $f(x_1,...,x_n,x_1^{\ast},...,x_n^{\ast})^*
f(x_1,...,x_n,x_1^{\ast},...,x_n^{\ast})$, one can assume that for
any n-tuple $T_i$ of operators  (or elements of any $C^*$-algebra)
$f(T_1,...,T_n,T_1^{\ast},...,T_n^{\ast})$ is a self-adjoint
operator (element).

It is sufficient to prove that there is $C>0$ such that for any
$\varepsilon
> 0$ and $K > 0$ there is a polynomial relation
\begin{equation}\label{soot2}
g(x_1,...,x_n, x_1^*, \ldots, x_n^*) = 0,
\end{equation}
$C$-boundedly representable in finite-dimensional $C^*$-algebra
and such that
\begin{equation}\label{razn}
\|f(T_1,...,T_n, T_1^*, \ldots, T_n^*) - g(T_1,...,T_n, T_1^*,
\ldots, T_n^*)\| < \varepsilon
\end{equation}
 for all operators $T_i$, $i=\overline{1, n}$, whose norms are not greater than  $K$.

Let   ${\mathcal A}$ be an  AF $C^*$-algebra in which (\ref{sootn'})
has a representation and $a^0_i\in {\mathcal A}$ be its elements
that satisfy (\ref{sootn'}). Then,  for any  $\delta$, there is a
finite-dimensional $\ast$-subalgebra
 $\mathcal B\subset \mathcal A$ such that $dist(a_i^0, \mathcal B) \le \delta < 1$, $i=\overline{1, n}$. Choosing $\delta$
small enough one can find such elements $b_n\in {\mathcal B}$ that
$\|f(b_1,...,b_n, b_1^*, \ldots, b_n^*)\|< \varepsilon/2$. Let $D
= \sup \|f(T_1,...,T_n, T_1^*, \ldots, T_n^*)\|$, where supremum
is taken over the set of all n-tuples of operators whose norms are
not greater than $K$.

Since ${\mathcal B}$ is finite-dimensional it has the form
${\mathcal B} = M_{k_1}\bigoplus\ldots \bigoplus M_{k_j}$ for some
$k_1, \ldots, k_j$, and hence we can think of its elements as of
matrices in $N$-dimensional space, where $N = k_1+\ldots + k_j$.
  Let $a = f(b_1,
..., b_n, b_1^*, \ldots, b_n^*)$,  $s_1 \ldots, s_N$ be the
eigenvalues of $a$ (repeated if multiple). By Lemma
\ref{mnogochlen} there exists a polynomial $h$ of one variable
such that $h(s_i)=0$, $i=\overline{1,N}$,  $|h(t) -
t|<\varepsilon$ when $|t|<D$. Since $a$ is self-adjoint we have
$h(a)=0$. Set
$$g(x_1,...,x_n, x_1^*, \ldots, x_n^*) = h(f(x_1,...,x_n, x_1^*,
\ldots, x_n^*)).$$ Then $g(x_1,...,x_n, x_1^*, \ldots, x_n^*)$ is
an element of free $\ast$-algebra on generators $x_1, \ldots, x_n$
and hence (\ref{soot2}) is a polynomial relation. We have $
g(b_1,...,b_n, b_1^*, \ldots, b_n^*) = h(a) = 0$, that follows
that (\ref{soot2}) is representable in finite-dimensional
$C^*$-algebra. Now, for any n-tuple $T_1,...,T_n$ of operators
with norms not greater than $K$, we have:
\begin{multline*}\|g(T_1, \ldots, T_n, T_1^*, \ldots, T_n^*) - f(T_1, \ldots, T_n, T_1^*, \ldots, T_n^*)\|
= \\ \|h(f(T_1, \ldots, T_n, T_1^*, \ldots, T_n^*)) - f(T_1,
\ldots, T_n, T_1^*, \ldots, T_n^*)\| \le \\sup_{|\lambda|\le
\|f(T_1, \ldots, T_n, T_1^*, \ldots, T_n^*)\|} |h(\lambda) -
\lambda| \le sup_{|\lambda|\le D}|h(\lambda) - \lambda| \le
\varepsilon.\end{multline*} It remains to note that (\ref{soot2})
is $(C+1)$-representable in a finite-dimensional space, where $C =
max_{i\in \{1, \ldots, n\}} \|a_i^0\|$. \qed

\medskip

It would be interesting to know if the inverse assertion  is true.

\begin{question}\label{question} Suppose that a polynomial relation is a limit of
polynomial relations  $C$-boundedly representable in
finite-dimensional $C^*$-algebras. Is it true that this relation
is representable in an AF-algebra?
\end{question}

\section{Universal $C^*$-algebras of relations (\ref{osnsootn})}


\subsection{$\Sigma_n(A)$ for type I $C^*$-algebras}

Recall that
a $C^*$-algebra $A$ is called a {\it CCR-algebra} if for any its
non-zero irreducible representation $(H, \pi)$ the set  $\pi(A)$
coincides with the set $K(H)$ of all compact operators on $H$.

We will use the definition of type I $C^*$-algebra which is given
in terms of composition series. The following two definitions can
be found in \cite{Sakai} or \cite{Glimm}.

\begin{opred} Let $A$ be a $C^*$-algebra. An increasing family of closed two-sided ideals
 $(I_{\rho})_{0\le \rho \le\alpha}$
of $A$, indexed by the ordinals between 0 and a certain ordinal
$\alpha$, is called a composition series of $A$ if it satisfies the
following conditions: \item 1) $I_0 = 0, I_{\alpha}= A$; \item 2) if
$\rho \le \alpha$ is a limit ordinal, then $I_{\rho} =
\overline{\bigcup_{\rho'<\rho}I_{\rho'}}$.
\end{opred}

\begin{opred}\label{2opr} A $C^*$-algebra $A$ is called type I $C^*$-algebra
if it has a composition series $(I_{\rho})_{0\le \rho \le\alpha}$
such that all $I_{\rho+1}/I_{\rho}$ are CCR-algberas.
\end{opred}

We couldn't find a reference for the following result.
\begin{prop}\label{reg}
Any (unital) type I $C^*$-algebra has a finite-dimensional
representation.
\end{prop}
{\bf Proof.} Let $A$ be a unital type I $C^*$-algebra and
$(I_{\rho})_{0\le \rho \le\alpha}$ be its composition series.
Suppose that $\alpha$ is a limit ordinal and hence $A =
\overline{\bigcup_{\rho<\alpha}I_{\rho}}$. Then the unit of $A$ is a
limit of sequence of elements that belong to ideals of $A$. It is
impossible because the set of invertible elements of $A$ is open.
Thus  $\alpha$ is not a limit ordinal and hence the composition
series has an ideal $I_{\alpha-1}$. Then $A/I_{\alpha-1}$ is a
unital CCR-algebra. Hence all its irreducible representations are
finite-dimensional (because the image of the unit should be
compact). Let $\rho$ be any irreducible representation of
$A/I_{\alpha-1}$, $q: A\to A/I_{\alpha-1}$ be the canonical
epimorphism. Then the composition $\rho{\circ}q$ gives a
finite-dimensional representation of $A$.
 \qed

Recall that, for any $C^*$-algebra $A$, we denote by $\Sigma_n(A)$
the set of those  $\lambda$ for which {\it in this algebra} there
exist  $n$  projections whose sum is $\lambda 1$.

\begin{theorem}\label{4Pr} Let $A$ be a type I $C^*$-algebra. Then  $\Sigma_n(A)$ consists of finitely many rational numbers.
\end{theorem}

{\bf Proof.} By Proposition \ref{reg}, there exists a representation
$\pi$ of $A$ in a finite-dimensional space $H$. Denote by  $m$ its
dimension. Let $\lambda \in \Sigma_n(A)$.  Then there exists
projections $p_1, \ldots, p_n$ in $A$ such that $p_1+\ldots + p_n =
\lambda 1$.
  Calculating traces of left-hand side
and right-hand side of the equality $$\pi(p_1)+\ldots +\pi(p_n) =
\lambda \pi(1),$$ we get
$$tr (p_1)+\ldots + tr (p_n) =
\lambda m.$$ Since the trace of any projection in m-dimensional
space is a natural number not greater than m, we obtain that
$\lambda$ belongs to the finite set of rational numbers
$\{\frac{k_1+\ldots+k_n}{m}: k_i\le m, k_i\in \mathbb N, i=1,
\ldots,n\}$.  Since it is true for any $\lambda \in \Sigma_n(A)$, we
are done. \qed

\begin{cor}\label{netipa1} $\mathcal{P}_{n, \alpha_n}$ and $\mathcal{P}_{n, \beta_n}$ are not type I $C^*$-algebras.
\end{cor}
{\bf Proof.} Clearly $\lambda\in \Sigma_n$ implies $\lambda\in
\Sigma_n(\mathcal{P}_{n, \lambda})$.  Hence $\alpha_n\in
\Sigma_n(\mathcal{P}_{n, \alpha_n})$, $\beta_n\in
\Sigma_n(\mathcal{P}_{n, \beta_n})$. Since  $\alpha_n$ and $\beta_n$
are irrational  we obtain, by Theorem \ref{4Pr}, that
$\mathcal{P}_{n, \alpha_n}$ and $\mathcal{P}_{n, \beta_n}$ are not
type I $C^*$-algebras. \qed

  \begin{cor} The set of polynomial relations, representable in type I $C^*$-algebras,
  is not closed.
\end{cor}
{\bf Proof.} Let $\lambda_k\in \Lambda_n^1$ (definition of
$\Lambda_n^1$ is given in Introduction), $\lambda_k \to \alpha_n$.
Then the relations
\begin{equation}\label{diskr}
\{\sum_{i=1}^n p_i = \lambda_k 1, p_i^2 = p_i, p_i^* = p_i\}
\end{equation}
are representable in type I $C^*$-algebras $\mathcal{P}_{n,
\lambda_n}$ respectively. On the other hand the relation
$\{\sum_{i=1}^n p_i = \alpha_n, p_i^2 = p_i, p_i^* = p_i\}$ is a
limit of the relations (\ref{diskr}), but, since $\alpha_n$ is
irrational, it cannot be represented in a type I $C^*$-algebra by
Theorem \ref{4Pr}. \qed

\subsection{ $\Sigma_n(A)$ for $C^*$-algebras with a trace, and
classification of $P_{n, \lambda}$, $\lambda\in \Sigma_n$}

\begin{theorem}\label{trace}
All $C^*$-algebras  $\mathcal{P}_{n, \lambda}$, where $\lambda\in
\Sigma_n$, have a trace.
\end{theorem}
{\bf Proof.} For any rational $\lambda\in \Sigma_n$, the relation
(\ref{osnsootn}) has a finite-dimensional representation
(\cite{Yu.S.}) and hence a representation in a $C^*$-algebra with
a trace.  Moreover, for any $\lambda\in \Sigma_n$,   each
representation of the relation (\ref{osnsootn}) is 1-bounded.
Since any irrational $\lambda \in \Sigma_n$ is a limit of rational
numbers from $\Sigma_n$, the relation (\ref{osnsootn}) belongs to
the closure of the set of polynomial relations $1$-representable
in  $C^*$-algebras with a trace.  Hence, by Theorem \ref{pr-sled}
it is representable in a $C^*$-algebra with a trace.

Since any representation of the relation (\ref{osnsootn}) defines
a representation of $C^*$-algebra $\mathcal{P}_{n, \lambda}$ we
get that for any $\lambda \in \Sigma_n$, $\pnl$ has a
$\ast$-homomorphism $\pi$ to a $C^*$-algebra $A$ with a trace. Let
us denote this trace by $\tau$. Setting $\tau_1(a) = \tau(\pi(a))$
for any $a \in \pnl$, we get a trace on $\mathcal{P}_{n,
\lambda}$. It is non-zero because $\tau_1(1_{\pnl}) = 1$. \qed

\begin{theorem}\label{2Pr} Let $A$ be a separable $C^*$-algebra with a trace. Then the set
 $\Sigma_n(A)$ is countable.
\end{theorem}

{\bf Proof.} The set $P(A)$ of all projections in $A$ is a subset
of separable metric space and hence is separable. Let $\tilde
P=\{p_k, k\in \mathbb N\}$ be a countable dense subset of $P(A)$.
Let $\lambda\in \Sigma_n(A)$. Then there exist  $q_i \in P(A),
i=\overline{1, n},$ such that $\sum_{i=1}^n q_i = \lambda 1$.
Since $\tilde P$ is dense in $P(A)$ we can find such $p_{k(i)}\in
\tilde P$ that
 $\|p_{k(i)}-q_i\| < 1$. This implies that the projections
$p_{k(i)}$ and $q_i$ are equivalent and hence $\tau(p_{k(i)}) =
\tau(q_i)$, where $\tau$ is a trace on $A$. It follows that
$\lambda = \sum_{i=1}^n \tau(p_{k(i)})$. But the set of all
n-tuples of projections from $\tilde P$ is countable. Hence
$\Sigma_n(A)$ is countable. \qed

 \medskip

 \begin{remark} In the absence of a trace the theorem is not true even for separable simple nuclear
 $C^*$-algebra $A$.  As an example one can take $O_2$ (see \cite{Yu.S.}).
\end{remark}

Consider now the problem of classification of the family $\pnl$,
$\lambda \in \Sigma_n$. We don't know in general when $C^*$-algebras
$\pnl$ and $P_{n, \mu}$ are isomorphic. It is natural to conjecture
that it happens only when $\mu = \lambda$ or $\mu = n-\lambda$. For
$n<5$ it is true. For $n\ge 5$,
 the invariant $\Sigma_n(A)$ helps to prove  that
among these $C^*$-algebras there is continuum of pairwise
nonisomorphic ones. Even more strongly

\begin{theorem}\label{neisom}
Let $E\subset \Sigma_n$ have the cardinality of continuum. Then
among $C^*$-algebras $\mathcal{P}_{n, \lambda}, \lambda\in E,$
there is continuum of pairwise nonisomorphic ones.
\end{theorem}
{\bf Proof.}  Represent the set of all $C^*$-algebras
  $\mathcal{P}_{n, \lambda}$, $\lambda\in E,$ as the union of classes of pairwise isomorphic $C^*$-algebras.
  Let $\{K_i: i\in I\}$ be the set of all these classes.  For all $C^*$-algebras $A$ from one equivalence class  $K_i$
the set $\Sigma_n(A)$ is the same, so we can denote it by
$\Sigma_n(K_i)$.
 By Theorems \ref{trace} and
\ref{2Pr},  for any $\lambda\in \Sigma_n$, the set
$\Sigma_n(\mathcal{P}_{n,
 \lambda})$ is countable and hence $\Sigma_n(K_i)$ is countable for any $i\in I$.
Since clearly $\lambda \in \Sigma_n(\mathcal{P}_{n, \lambda})$ we
have $E = \cup_{i\in I}\Sigma_n(K_i)$. Thus $\cup_{i\in
I}\Sigma_n(K_i)$ has the cardinality of continuum and each
$\Sigma_n(K_i)$ is countable. It follows that  $I$ has the
cardinality of continuum. \qed

\subsection{$\Sigma_n(A)$ for UHF-algebras}
A $C^*$-algebra is called {\it uniformly hyperfinite} (UHF, in
short) if it is the union of an increasing  net of unital
subalgebras isomorphic to full matrix algebras.  For such
$C^*$-algebras the set $\Sigma_n(A)$ can be written explicitly.

\begin{theorem}\label{hyperf}
Let $A =  \overline{\bigcup A_i}$ be a UHF-algebra with $A_i\cong
M_{k_i}$. Then the set $\Sigma_n(A)$ consists of all numbers
$\lambda \in \Sigma_n$ of the form $p/q$, where $q|k_j$ for some
$j$.
\end{theorem}

{\bf Proof.} Let $p/q \in \Sigma_n$,  $q|k_j$ for some $j$. Clearly
one can assume that $p/q$ is irreducible fraction. By definition,
$(p/q)1$ is the sum of $n$ projections in $q$-dimensional space
(\cite{LinAlg}). If $q|k_j$ for some $j$ then there is an embedding
of  $M_q$ into $M_{k_j}$ and hence into $A$. Thus we have $n$
projections in $A$ with sum $(p/q)1$.

Prove that no other  numbers can belong to $\Sigma_n(A)$. Clearly
$\Sigma_n(A) \subset \Sigma_n$ because $A$ can be enclosed in
$B(H)$. Suppose $\sum_{i=1}^n r_i = \lambda 1$, where all $r_i$
are projections. Since $\bigcup_j A_j$ is dense in $A$ there are
$j\in \mathbb N$ and elements $a_i \in A_j$, $i=\overline{1, n}$,
such that $\|a_i-r_i\| < 1$, $i=\overline{1, n}$. Using a standard
trick with functional calculus (see, for example, \cite{Davidson},
section III.3) one can change all $a_i$ by projections $d_i$.
  Then $d_i$ is
equivalent to $r_i$, $i=\overline{1, n}$,  their (normalized)
traces are the same and we get
 $$\lambda = \sum_{i=1}^n \tau (r_i) = \frac{\sum_{i=1}^n tr(d_i)}{k_j}.$$
 \qed

In terms of supernatural numbers (for definition, see, for
example, \cite{Davidson}) the theorem can be reformulated in the
following way. Let us say that $n$ divides supernatural number
$\delta(A) = \prod p^{\varepsilon_p}$ if the exponent of every
prime divisor $q$ of $n$ in the factorization of $n$  is not
greater than $\varepsilon_q$. Then Theorem \ref{hyperf} says that
$\Sigma_n(A)$ consists of all rational numbers  $\lambda \in
\Sigma_n$ whose denominators divide the supernatural number of
$A$.

As  we know when $A$ is $B(H)$ or type I $C^*$-algebra the set
$\Sigma_n(A)$ is closed. It follows from Theorem \ref{hyperf} that
for UHF-algebras it is not true (this answers a question of the
authors of \cite{Yu.S.}).

\subsection{Nuclearity and exactness}

$C^*$-algebra $A$ is {\it nuclear} if, for any $C^*$-algebra $B$,
there is only one $C^*$-norm on the algebraic tensor product
$A{\odot}B$. For the theory of representations the most important
characterization of this class of $C^*$-algebras is the following:
a $C^*$-algebra is nuclear if and only if any its
factor-representation generates hyperfinite factor.

Our aim in this subsection is to prove that for large $n$,
$\mathcal{P}_{n, \lambda}$ is non-nuclear for the "most" of points
$\lambda \in (\alpha_n;\beta_n)$.
 Moreover we will show that for $n>10$, $(\alpha_n;\beta_n)$ contains a
 subinterval
 $I_n$ such that for any $\lambda\in I_n$,
  $\mathcal{P}_{n, \lambda}$ doesn't belong to much larger class of exact $C^*$-algebras.

 A $C^*$-algebra $A$ is called {\it exact} if, for any short exact sequence $$0\to J\to B\to C\to 0,
$$ the sequence $$0\to A{\otimes}J\to A{\otimes}B\to
A{\otimes}C\to 0 $$ is also exact. By  $\otimes$ we denote the
minimal tensor product.

It is well known that the class of all nuclear $C^*$-algebras is
contained in the class of all exact $C^*$-algebras. Recall also that
both classes are closed under taking ideals and quotients and that
the class of exact algebras is closed under taking  closed
$\ast$-subalgebras. All this information can be found in
\cite{Kirch}, \cite{Kirch2}.

Below $C^{\ast}(\mathbb{F}_2)$ means the group $C^*$-algebra of
the  free group on two generators.

\begin{lemma}\label{trikarty}
In infinite-dimensional Hilbert space there exist 3 projections $P,
Q, R$ generating non-exact  $C^*$-algebra.
\end{lemma}

{\bf Proof.} Let $A$ be the universal $C^*$-algebra generated by 3
projections $p_1, p_2, p_3$ without any relations. From
[\cite{Ostr-Sam}, Theorem 54, Proposition 66] it follows that
there exists a closed ideal  $J$ of $A$ such that $A/J\cong
M_n{\otimes}C^{\ast}(\mathbb{F}_2)$ for some $n\in \mathbb N
\bigcup \{\infty\}$ (here $M_{\infty}$ means the algebra of all
compact operators).

This implies that  $A$ is non-exact. Indeed if $A$ is exact then
any its quotient is exact. On the other hand
$C^{\ast}(\mathbb{F}_2)$ is non-exact (\cite{Wass}) and hence
$M_n{\otimes}C^{\ast}(\mathbb{F}_2)$ is non-exact because it
contains non-exact $C^*$-algebra  $C^{\ast}(\mathbb{F}_2)$ as a
closed $\ast$-subalgebra.

Now let $\pi$ be the universal representation of $A$. Set $P =
\pi(p_1), Q=\pi(p_2), R=\pi(p_3)$, then $C^*$-algebra generated by
them is isomorphic to $A$ and hence is non-exact. \qed

\begin{theorem}\label{6Pr} For each $n>6$, there exists a nonempty subset
$I_n\in \Sigma_n$ such that for any $\lambda\in I_n$, the
$C^*$-algebra $\mathcal{P}_{n, \lambda}$ is not exact.

If $n>10$ then $I_n \supset [5;n-5]$.
\end{theorem}

{\bf Proof.} Consider such $\lambda\in \Sigma_n$ that   $\lambda -3
\in \Sigma_{n-6}$. The set $\Sigma_n\cap(\Sigma_{n-6} + 3)$ of all
such points we denote by $E_n$. Using  (\ref{KontsyOtrezka}) it is
easy to check that $E_n\neq \emptyset$, for $n> 6$, and that
$\beta_{n-6}+3 < \beta_n$, $\alpha_{n-6}+3> \alpha_n$ for any $n
> 10$, whence we get that  $E_n$ contains the closed interval $I_n =
[\alpha_{n-6}+3;\beta_{n-6}+3]$. Since $[2;n-2] \subset
[\alpha_n;\beta_n]$ we get $I_n\supset [5,n-5]$ for any $n>10$.

 Let $\lambda_1 \in \Sigma_{n-6}$
and let $\pi$ be arbitrary representation of the $C^*$-algebra
$\mathcal{P}_{n-6, \lambda_1}$. Define a representation $\tilde \pi$
of $\mathcal{P}_{n, \lambda}$, where $\lambda = \lambda_1 + 3\in
E_n$, in the following way. Set

\medskip

$\tilde \pi(p_k) = \pi(p_k)$ for $1\le k \le n-6$,

\medskip

 $\tilde
\pi(p_{n-5}) = P$, $\;\tilde \pi(p_{n-4}) = 1-P$,

\medskip

$\tilde \pi(p_{n-3}) = Q$, $\;\tilde \pi(p_{n-2}) = 1-Q$,

\medskip

$\tilde \pi(p_{n-1}) = R$, $\;\tilde \pi(p_{n}) = 1-R$,

\medskip

\noindent  where $P, Q, R$ are  projections constructed in  Lemma
\ref{trikarty}. Since $C^*$-algebra generated by them is a
subalgebra of $\tilde \pi (\mathcal{P}_{n, \lambda})$ it is
isomorphic to some subalgebra of the quotient
 $ \mathcal{P}_{n, \lambda}/Ker{\tilde \pi}$.
Hence $\mathcal{P}_{n, \lambda}$ is not exact because any quotient
of exact $C^*$-algebra is exact and any subalgebra of an exact
$C^*$-algebra is exact (\cite{Kirch2}). \qed

\begin{remark}
 Since $[\alpha_n;\beta_n] \subset [1;n-1]$ for any $n$, and
$I_n$ contains  $[5;n-5]$ for $n>10$, we can say that for large $n$,
$I_n$ contains "almost whole" $\;[\alpha_n;\beta_n]$.
\end{remark}

Now we are going to prove that the set of points $\lambda$ such that
$\mathcal{P}_{n, \lambda}$ is not nuclear, is strictly larger than
$I_n$.

Let us denote by $f$ the map from the interval $(\alpha_n;\beta_n)$
onto itself given by the formula  $f(\lambda) = n-1- 1/(\lambda-1)$.
Let $\mathcal{S}(f)$ be the group (with the composition as a group
multiplication) generated by $f$, that is the group of all (positive
and negative) powers of the map $f$.

 It was proved in \cite{Yu.S.}  that if  $\lambda_1 = f(\lambda_2)$
then the categories of representations of $C^*$-algebras
$\mathcal{P}_{n, \lambda_1}$ and $\mathcal{P}_{n, \lambda_2}$ are
equivalent. This means that there exist a bijection $\pi\to
\tilde{\pi}$ between the sets of representations of these
$C^*$-algebras  and, for any   $\pi_1,\pi_2\in \it{Rep}(P_{n,
\lambda_1})$, the linear bijection $F_{\pi_1,\pi_2}$ from the
space of intertwining operators $W(\pi_1,\pi_2)$ onto
$W(\tilde{\pi_1},\tilde{\pi_2})$. Moreover the map
$F_{\pi_1,\pi_2}$ and its inverse are continuous in WOT. Also, as
for any functor, if $T\in W(\pi_1,\pi_2)$, $S\in W(\pi_2,\pi_3)$
then $F_{\pi_1,\pi_3}(ST) = F_{\pi_2,\pi_3}(S)F_{\pi_1,\pi_2}(T)$.
It follows that if $\pi\in \it{Rep}(P_{n, \lambda_1})$ is a
factor-representation then $\tilde \pi\in \it{Rep}(P_{n,
\lambda_2})$ is also a factor-representation.

\begin{theorem}
Let $n>6$. $C^*$-algebra $\mathcal{P}_{n, \lambda}$ is not nuclear
for every $\lambda \in (\alpha_n;\beta_n)$ whose orbit of the
action of  $\mathcal{S}(f)$ intersects $I_n$.
\end{theorem}
{\bf Proof.} Let  $\lambda_1 = f(\lambda_2)$.  We have to prove that
if one of $C^*$-algebras  $\mathcal{P}_{n, \lambda_i}$, $i=1,2$, is
nuclear then the second one is also nuclear.

Suppose that $\mathcal{P}_{n, \lambda_1}$ is not nuclear. Then there
is its factor-representation $\pi$ which is not hyperfinite, that
means that the closure in WOT of $\pi(\mathcal{P}_{n, \lambda_1})$
is not hyperfinite. By Connes's theorem (\cite{Connes}) its
commutant $\pi(\mathcal{P}_{n, \lambda_1})'$ is not hyperfinite too.
But because of mentioned above $\pi(\mathcal{P}_{n, \lambda_1})' =
W(\pi,\pi)$ is isomorphic, as $W^*$-algebra, to
$\tilde{\pi}(\mathcal{P}_{n, \lambda_2})' =
W(\tilde{\pi},\tilde{\pi})$. Hence $\tilde{\pi}(\mathcal{P}_{n,
\lambda_2})'$ is not hyperfinite. Applying again Connes's theorem we
obtain that  $\mathcal{P}_{n, \lambda_2}$ is not nuclear because it
has a factor-representation which is not hyperfinite.

Thus all $C^*$-algebras  $\mathcal{P}_{n, \lambda}$, with $\lambda$
from one orbit of the action of $\mathcal{S}(f)$, are nuclear or
non-nuclear simultaneously. Now it remains to apply Theorem
\ref{6Pr}. \qed

\subsection{Concluding remarks} We will mention some additional results and questions about  $\mathcal{P}_{n,
\lambda}$.

 1) Stability.

\noindent  Let $\delta
> 0$. An n-tuple of operators $T_1,...,T_n$ is called
a $\delta$-representation of the relation (\ref{sootn}) if
$$ \|f(T_1,...,T_n)\|\le \delta.$$
 The relation (\ref{sootn})
is called {\it stable} (see \cite{Loring}) if for any $\varepsilon
> 0$ there is  $\delta > 0$ such that if  $T_1,...,T_n$
is a $\delta$-representation of this relation then there exists its
representation  $S_1,...,S_n$ in the same space that satisfies the
condition $\|T_i-S_i\| < \varepsilon$.

\begin{theorem}\label{ustoj}
For any $\lambda \in [\alpha_n, \beta_n]$, the relation
(\ref{osnsootn}) is not stable.
\end{theorem}
{\bf Proof.} It suffices to show that for each $\delta
> 0$ there exists a $\delta$-representation of the
relation (\ref{osnsootn}) in a Hilbert space $H$ but there  are no
representations of (\ref{osnsootn}) in $H$.  We will consider
separately the case when $\lambda$ is rational and the case when
it is irrational.

Let $\lambda\in [\alpha_n, \beta_n]$ be irrational.   Let
$\lambda'\in \Sigma_n$ be  rational  and
 $|\lambda' - \lambda| < \delta$.
By \cite{Yu.S.}, there exist  projections $P_1,...,P_n$ in a
finite-dimensional Hilbert space $H$ such that  $\sum_{i=1}^nP_i =
\lambda' 1$. Clearly they define a $\delta$-representation of
(\ref{osnsootn}). Suppose that (\ref{osnsootn}) has a
finite-dimensional representation. Then there exist projections
 $Q_1, \ldots, Q_n\in B(H)$ such that
$Q_1+\ldots + Q_n = \lambda 1$. Taking a trace in both sides of
this equality we get that  $\lambda$ is rational. Hence
(\ref{osnsootn}) doesn't have any representation in $H$ and we are
done.

Now let $\lambda \in [\alpha_n, \beta_n]$ be rational.  Let $\lambda
= \frac{p}{q}$, where $\frac{p}{q}$ is irreducible fraction. There
exists a rational number $\lambda'$ of the form $\lambda' =
\frac{r}{p^{m}}$ such that $|\lambda' - \lambda|\le \delta$. By
\cite{Yu.S.} there exist projections $P_1, \ldots, P_n$ in $\mathbb
C^{p^{m}}$ such that $P_1+\ldots P_n = \lambda' 1$. Hence $P_1,
\ldots, P_n$ define a $\delta$-representation of (\ref{osnsootn}).
Suppose that this relation has some representation $Q_1, \ldots,
Q_n$ in $\mathbb C^{p^{m}}$. Then $tr Q_1+ \ldots tr Q_n = \lambda
p^{m}$ whence we get that $\lambda p^{m}\in \mathbb{Z}$. It follows
that $q$ and $p$ are not coprimes - in  contradiction to the
assumption. Hence (\ref{osnsootn}) doesn't have any representation
in $\mathbb C^{p^{m}}$ and we are done.

 \qed

 2) Simplicity.

\noindent For "most" $\lambda\in \Sigma_n$, it is easy to prove that
$\pnl$ is not simple. Namely

(i) if $\lambda\in \Lambda_n^i$, $i=1, 2$, then $\pnl$ is full
matrix algebra or sum of full matrix algebras (\cite{Yu.S.},
Theorem 4).

(ii) if $\lambda\in [\alpha_n, \beta_n]$ is rational then
$\mathcal{P}_{n, \lambda}$ is not simple because it is not
finite-dimensional (even not type I) but has a finite-dimensional
representation (\cite{LinAlg}).

(iii) if $\lambda\in \Sigma_{n-1}\bigcap \Sigma_n$ then
$\mathcal{P}_{n, \lambda}$ is not simple. Indeed we can define a
$\ast$-homomorphism $\pi : \mathcal{P}_{n, \lambda}\to
\mathcal{P}_{n-1, \lambda}$ setting $\pi(p_i)=q_i$,
$i=\overline{1, n-1}$, $\pi(p_n)=0$, where $p_1,\ldots, p_n$ and
$q_1, \ldots, q_{n-1}$ are generators of
 $\mathcal{P}_{n, \lambda}$ and
 $\mathcal{P}_{n-1, \lambda}$ respectively, and take its kernel.

 So the question if $\pnl$ is simple remains open for
irrational numbers from $ [\alpha_n,
\beta_n]\backslash[\alpha_{n-1}, \beta_{n-1}]$.


\medskip

 3) K-theory.

 \noindent  D. Hadwin (private communication) proved that for
any $\lambda \in [\alpha_n, \beta_n]$ the group $K_0(\pnl)$ contains
$\mathbb Z^n$ as a direct summand.

\medskip

 4) It would be interesting to calculate $\Sigma_n(A)$
for any von Neumann algebra $A$. It is not difficult to see that
the problem can be reduced to the case when $A$ is a factor and
for factors this problem is not trivial only in the case when $A$
is $II_1$-factor. The most intriguing is the case when $A$ is the
 hyperfinite $II_1$-factor. If Question
\ref{question} (Chapter 1) had positive answer it would be easy to
prove that for hyperfinite $II_1$-factor (and therefore for any
infinite-dimensional factor) $\Sigma_n(A) = \Sigma_n$.

\end{document}